\documentclass[titlepage,11pt]{article}

\usepackage[margin=1.2in]{geometry}
\usepackage{latexsym}
\usepackage{amsfonts}
\usepackage[leqno]{amsmath}

\usepackage{thmtools}
\usepackage{thm-restate}

\makeatletter
\newcommand{\leqnomode}{\tagsleft@true}
\newcommand{\reqnomode}{\tagsleft@false}
\makeatother

\newcommand{\dichi}{\vec{\chi}}

\usepackage{tikz}
\usepackage{float}
\usepackage{url}
\usetikzlibrary{arrows}
\usetikzlibrary{decorations.markings}
\usepackage[symbol]{footmisc}

\usepackage{xspace}
\usepackage[color=Olive]{todonotes}

\usepackage[colorlinks=true,allcolors=blue]{hyperref}
\usepackage{cleveref}

\newcommand\blackslug{\hbox{\hskip 1pt \vrule width 4pt height 8pt depth 1.5pt
        \hskip 1pt}}
\newcommand\bbox{\hfill \quad \blackslug \bigbreak}

\usepackage{authblk}

\title{A counterexample to a conjecture about triangle-free induced subgraphs of graphs with large chromatic number}

\author[1]{Alvaro Carbonero}
\author[1]{Patrick Hompe}
\author[2]{Benjamin Moore}
\author[1]{Sophie Spirkl\thanks{Emails: (ar2carbonerogonzales, phompe, sspirkl)@uwaterloo.ca, brmoore@iuuk.mff.cuni.cz \\
We acknowledge the support of the Natural Sciences and Engineering Research Council of Canada (NSERC), [funding reference number RGPIN-2020-03912]. Cette recherche a été financée par le Conseil de recherches en sciences naturelles et en génie du Canada (CRSNG), [numéro de référence RGPIN-2020-03912].
Benjamin Moore is supported by the ERC-CZ project LL2005 (Algorithms and complexity within and
beyond bounded expansion) of the Ministry of Education of Czech Republic.
}}

\affil[1]{University of Waterloo, Department of Combinatorics and Optimization, Waterloo, Canada}

\affil[2]{Charles University, Institute of Computer Science, Prague, Czech Republic}

\date{\today}

\newtheorem{thm}{Theorem}[section]

\newtheorem{conjecture}[thm]{Conjecture}
\newtheorem{lemma}[thm]{Lemma}

\newtheorem{question}[thm]{Question}

\newcommand{\Proof}{\noindent{\bf Proof.}\ \ }

\begin{document}
\maketitle
\begin{abstract}
We prove that for every $n$, there is a graph $G$ with $\chi(G) \geq n$ and $\omega(G) \leq 3$ such that every induced subgraph $H$ of $G$ with $\omega(H) \leq 2$ satisfies $\chi(H) \leq 4$.

This disproves a well-known conjecture. Our construction is a digraph with bounded clique number, large dichromatic number, and no induced directed cycles of odd length at least 5. 
\end{abstract}

\section{Introduction and preliminaries}

In this paper, we disprove the following conjecture (its origin appears somewhat unclear;\footnote{Recently, Vojt\v{e}ch R\"odl pointed out to us that the problem had first appeared in \cite{nesetril} attributed to Fred Galvin and Vojt\v{e}ch R\"odl.} it is attributed to Louis Esperet in \cite{conjecture2}, while the authors of \cite{conjecture1} state that ``we
could not find a reference''):

\begin{conjecture}[\cite{conjecture2, conjecture1}] \label{conj:main}
For all $k, r \in \mathbb{N}$ there is an $n \in \mathbb{N}$ such that for every graph $G$ with $\chi(G) \geq n$ and $\omega(G) \leq k$, there is an induced subgraph $H$ of $G$ with $\chi(H) \geq r$ and $\omega(H) = 2$.  
\end{conjecture}

Here, $\chi(G)$ denotes the \emph{chromatic number} of a graph $G$ and $\omega(G)$ denotes the \emph{clique number}. This conjecture is the induced-subgraph analogue of the following theorem: 

\begin{thm}[\cite{rodl}]
For every $r \in \mathbb{N}$ there is an $n \in \mathbb{N}$ such that for every graph $G$ with $\chi(G) \geq n$, there is a (not necessarily induced) subgraph $H$ of $G$ with $\chi(H) \geq r$ and $\omega(H) = 2$.  
\end{thm}

We will show: 

\begin{restatable}{thm}{maina}
\label{thm:main1}
For every $n \in \mathbb{N}$, there is a graph $G$ with $\chi(G) \geq n$ and $\omega(G) \leq 3$ such that  every induced subgraph $H$ of $G$ with $\omega(H) \leq 2$ satisfies $\chi(H) \leq 4$.
\end{restatable}

This is answers Conjecture \ref{conj:main} in the negative for all $r \geq 5$. The following shows that the case $r = 4$ is the only case of Conjecture \ref{conj:main} which remains open:
\begin{thm}[\cite{oddholes}]  \label{thm:oddholes}
There is a function $f : \mathbb{N} \rightarrow \mathbb{N}$ such that for every graph $G$ with no induced cycle of odd length at least 5, we have $\chi(G) \leq f(\omega(G)).$
\end{thm}
Letting $f$ as in Theorem \ref{thm:oddholes}, it follows that every graph $G$ with $\chi(G) > f(\omega(G))$ contains an induced cycle of odd length at least 5, and therefore contains an induced subgraph $H$ with $\omega(H) = 2$ and $\chi(H) = 3$. 

Our construction is based on a construction of \cite{kt} and produces a digraph with large dichromatic number, which we define below. Throughout this paper, we only consider simple digraphs $D$, that is, for every two distinct vertices $u$ and $v$, the digraph $D$ contains either an edge from $u$ to $v$, or an edge from $v$ to $u$, or neither; but not both. For a digraph, we write $uv$ for an edge from $u$ to $v$. 
Given a digraph $D$, we define its \emph{underlying undirected graph $G$} to be that graph with $V(G) = V(D)$ and in which $u, v \in V(G)$ are adjacent if $D$ contains an edge from $u$ to $v$ or from $v$ to $u$. The \emph{clique number $\omega(D)$} of a digraph $D$ is defined as the clique number of the underlying undirected graph of $D$.

An analogue of chromatic number for directed graphs was introduced in \cite{dichi1, dichi2}. A digraph is \emph{acyclic} if it contains no directed cycle. For $k \in \mathbb{N}$, a \emph{$k$-dicoloring} of a digraph $D$ is a function $f: V(D) \rightarrow \{1, \dots, k\}$ such that for every $i \in \{1, \dots, k\}$, the induced subdigraph of $D$ with vertex set $\{v \in V(D): f(v) = i\}$ is acyclic. The \emph{dichromatic number $\dichi(D)$} is the smallest integer $k$ such that $D$ has a $k$-dicoloring. 

We show that the digraph analogue of Theorem \ref{thm:oddholes} does not hold: 
\begin{restatable}{thm}{mainb}
\label{thm:main2}
For every $n$, there is a digraph $D$ with $\dichi(D) \geq n$, $\omega(D) \leq 3$ and with no induced directed cycle of odd length at least 5. 
\end{restatable}

\section{The construction}
We construct a sequence of digraphs $\{D_n\}$ as follows. Let  $D_1$ be the digraph with a single vertex. For $n \ge 2$, we take $n-1$ disjoint copies of the digraph $D_{n-1}$ and call them $D_{n-1}^1, \dots, D_{n-1}^{n-1}$. Let $\mathcal{T}$ be the set of all sequences $T = (x_1, \dots, x_{n-1})$ with $x_i \in V(D_{n-1}^{i})$ for all $i \in \{1, \dots, n-1\}$. Now, for every $T = (x_1, \dots, x_{n-1}) \in \mathcal{T}$ we create a vertex $v_T$ and for every $i \in \{1, \dots, n-1\}$, we add an edge from $x_i$ to $v_T$. The resulting digraph with vertex set $$V(D_{n-1}^1) \cup \dots \cup V(D_{n-1}^{n-1}) \cup \{v_T : T \in \mathcal{T}\}$$ and edge set $$E(D_{n-1}^1) \cup \dots \cup E(D_{n-1}^{n-1}) \cup \{x_iv_T : i \in \{1, \dots, n-1\}, T = (x_1, \dots, x_{n-1}) \in \mathcal{T}\}$$ is called $D_n$.  

We note that the graph $D_n$ is the graph of red edges in the proof of Theorem 3 of \cite{kt}, where the following was proved: 
\begin{lemma}[\cite{kt}] \label{thm:red}
For all $n \in \mathbb{N}$, we have: 
\begin{itemize}
    \item $D_n$ is acyclic; 
    \item for every two vertices $u, v \in V(D_n)$ there is at most one directed path from $u$ to $v$ in $D_n$. 
\end{itemize}
\end{lemma}
\Proof
We include a proof for completeness. For $n \geq 1$, let us define a partition of $V(D_n)$ into sets $T^n_1, \dots, T^n_n$ as follows: For $n=1$, let $T^1_1 = V(D_1)$. For $n > 1$ and $i \in \{1, \dots, n-1\}$, let $T^n_i$ be the union of the sets $T^{n-1}_i$ in $D_{n-1}^1, \dots, D_{n-1}^{n-1}$, and let $T^n_n$ be the set of remaining vertices (and thus $T^n_n$ is the set of vertices $v_T$ added when constructing $D_n$). 

By construction we have that for all $i \in \{1, \dots, n\}$, the set $T^n_i$ is a stable set and the only edges between $T^n_i$ and $T^n_1 \cup \cdots \cup T^n_{i-1}$ are edges from $T^n_1 \cup \cdots \cup T^n_{i-1}$ to $T^n_i$. It follows that $D_n$ is acyclic, as desired.

For the second bullet, note that every edge is from $T^n_i$ to $T^n_{j}$ for some $i < j$. Now, suppose we have vertices $u,v$ such that there exists a directed path $P$ from $u$ to $v$. Then it follows that $u \in T^n_i$ and $v \in T^n_j$ for $i < j$, and the vertex set of $P$ is contained in $T^n_i \cup \dots \cup T^n_j$. Let $H$ be the copy of $D_{j-1}$ that $u$ is contained in from the construction of $D_j$. By construction, every edge of $D_n$ with one end in $H$ and one end $x$ in $V(D_n) \setminus V(H)$ satisfies $x \in T^n_k$ for some $k \geq j$. Since $v$ is the only vertex of $P$ in $T^n_j \cup T^n_{j+1} \cup \dots \cup T^n_n$, it follows that all vertices of $P \setminus v$ are contained in $H$. Note that $v$ has exactly one in-neighbor in $H$; let that in-neighbor be $w$. It follows that any directed path from $u$ to $v$ must go through $w$. By induction on $n$ (since $P \setminus v$ is contained in a copy of $D_{j-1}$ with $j \leq n$), we have that there is at most one directed path from $u$ to $w$, so it follows that there is at most one directed path from $u$ to $v$, as desired. This completes the proof.\bbox{}

We define the \emph{length} of a (directed) path as its number of edges. Now, we construct a sequence of digraphs $\{D_n'\}$ as follows. We take a copy of $D_n$, and create a new graph $D_n'$ with $V(D_n') = V(D_n)$, and the following edges. For every two vertices $u,v$ where there exists a directed path in $D_n$ from $u$ to $v$,
\begin{itemize}
    \item we add an edge from $u$ to $v$ if that path has length equal to $1$ modulo $3$; and 
    \item we add an edge from $v$ to $u$ if that path has length equal to $2$ modulo $3$.
\end{itemize}
From Lemma \ref{thm:red}, it follows that $D_n'$ is well-defined and a simple digraph. In our analysis, it will be useful to consider a partition of the edges of $D_n'$ into two sets, positive and negative, which we call the \emph{sign} of an edge. Let us call an edge \emph{positive} if it was added as a result of the first bullet above, and \emph{negative} if it was added as a result of the second bullet. Clearly, this is a partition of the edges of $D_n'$. Note that in particular, if $uv \in E(D_n)$, then the edge $uv$ is added to $D_n'$ according to the first bullet, and hence $D_n$ is a (non-induced) subdigraph of the positive edges of $D_n'$. 

\begin{lemma} \label{thm:2path}
Let $u, v, w \in V(D_n')$. If $uv$ and $vw$ are edges of $D_n'$ of the same sign, then $wu$ is an edge of $D_n'$ of the opposite sign. 
\end{lemma}
\Proof
Suppose first that $uv$ and $vw$ are positive edges. Then by definition there exists a path $P_1$ from $u$ to $v$ in $D_n$ with length equal to $1$ modulo $3$, and a path $P_2$ from $v$ to $w$ in $D_n$ with length equal to $1$ modulo $3$. Then clearly $P_3 = P_2 \cup P_1$ is a directed walk from $u$ to $w$, and since $D_n$ is acyclic by Lemma \ref{thm:red}, it follows that $P_3$ is the unique directed path from $u$ to $w$. Then $P_3$ has length equal to $2$ modulo $3$, so it follows that $wu$ is a negative edge, as desired.

Suppose instead that $uv$ and $vw$ are negative edges. Then there exists a path $P_1$ from $v$ to $u$ and a path $P_2$ from $w$ to $v$ such that $P_1$ and $P_2$ both have length equal to $2$ modulo $3$. Then clearly $P_3 = P_2 \cup P_1$ is a directed walk from $w$ to $u$, and since $D_n$ is acyclic by Lemma \ref{thm:red}, it follows that $P_3$ is the unique path from $w$ to $u$. Then $P_3$ has length equal to $1$ modulo $3$ and it follows that $wu$ is a positive edge, as desired. This completes the proof.\bbox{}

\begin{lemma} \label{thm:triangle}
Let $u, v, w \in V(D_n')$. Then not all of $uv,vw, uw$ are edges of $D_n'$.
\end{lemma}
\Proof
We only consider the case when $uw$ is positive; the case when $uw$ is negative is analogous. It follows that there is a directed path $P_1$ from $u$ to $w$ of length congruent to $1$ modulo $3$. By Lemma \ref{thm:2path}, we may assume that $uv$ and $vw$ do not have the same sign. We consider two cases. 

If $uv$ is negative, then $vw$ is positive. It follows that there is a directed path $P_2$ from $v$ to $u$ of length congruent to $2$ modulo $3$. Now $P_3 = P_2 \cup P_1$ is a directed walk and since $D_n$ is acyclic by Lemma \ref{thm:red}, a directed path, from $v$ to $w$. But $P_3$ has length congruent to $0$ modulo $3$, and so from the construction of $D_n'$, it follows that $v$ and $w$ are not adjacent in either direction, a contradiction. 

Now $uv$ is positive, and $vw$ is negative. It follows that there is a directed path $P_2$ from $w$ to $v$ of length congruent to $2$ modulo $3$. Now $P_3 = P_1 \cup P_2$ is a directed walk and since $D_n$ is acyclic by Lemma \ref{thm:red}, a directed path from $u$ to $v$. But $P_3$ has length congruent to $0$ modulo $3$, and so from the construction of $D_n'$, it follows that $v$ and $u$ are not adjacent in either direction, a contradiction. \bbox{}

Now, we are ready to prove our main theorem, which we restate.
\maina*
\Proof
Let $\{G_n\}$ be the sequence of graphs such that $G_n$ is the underlying undirected graph of $D_n'$. Then we claim that taking $G = G_n$ will show the desired result.

Indeed, we first show that $\chi(G_n) \ge n$. Since $D_n$ is a subgraph of $D_n'$, it suffices to show, by induction, that the underlying undirected graph $H_n$ of $D_n$ has chromatic number at least $n$ (which was also shown in \cite{kt}, and follows from the fact that the $n$-th Zykov graph \cite{zykov} is a subgraph of $H_n$; here we give the short proof for completeness). The base case is trivial. By induction, we know that the underlying undirected graphs $H_{n-1}$ of the $n-1$ copies of $D_{n-1}$ that were used to build $D_n$ all have chromatic number at least $n-1$. So, if we take a coloring of $H_n$ with colors $\{1, \dots, n-1\}$, it follows that for every $i \in \{1, \dots, n-1\}$, there exists a vertex $x_i \in V(D_{n-1}^i)$ which receives color $i$. Then, letting $T = (x_1, \dots, x_{n-1})$, the corresponding vertex $v_T$ must receive a  color not in $\{1, \dots, n-1\}$, and it follows that the coloring uses at least $n$ colors. Thus, $\chi(H_n) \ge n$ for all $n \ge 1$, as claimed.

Next, we claim that $\omega(G_n) \le 3$. Suppose not; then $G_n$ contains a clique $K$ of size 4. Let $u \in K$ with at least two outneighbors in the digraph induced by $K$ in $D_n'$ (which is possible, since the average outdegree in this four-vertex digraph is 1.5), and let $v, w$ be two outneighbours of $u$ in $K$. By symmetry, we may assume that $vw$ is an edge of $D_n'$. But now $uv, vw, uw$ are all edges of $D_n'$, contrary to Lemma \ref{thm:triangle}.

Now, suppose that we have an induced subgraph $H$ of $G_n$ with $\omega(H) \le 2$. If we look at the corresponding induced subdigraph $H'$ of $D_n'$, it follows by Lemma \ref{thm:2path} that $H'$ does not contain a directed 2-edge path with both edges of the same sign as a subdigraph. Thus we can partition the vertices of $H'$ (and $H$) into two sets $A,B$ such that every vertex in $A$ is not the head of a positive edge and every vertex in $B$ is not the tail of a positive edge. Then note that there can be no positive edges between any two vertices in $A$, and also there are no positive edges between any two vertices in $B$. Likewise, we can find a similar partition $V(H') = A' \cup B'$ for the negative edges. Now $(A \cap A', A \cap B', B \cap A', B \cap B')$ is a partition of the vertices of $H$ into four stable sets, and thus $\chi(H) \le 4$, as claimed. This completes the proof.\bbox{}

The collection of digraphs $\{D_n'\}$ also gives the following result on $\dichi$-boundedness, which we restate.

\mainb*
\Proof
We claim that taking $D = D_{4n}'$ gives the desired result. Indeed, we know from the previous proof that $\omega(D) \le 3$. Furthermore, suppose that $D$ contains an induced odd directed cycle of length at least $5$. Then it follows that there exist two consecutive edges in that cycle of the same sign; but now Lemma \ref{thm:2path} gives a third edge which contradicts the fact that the cycle is induced.

It remains to show that $\dichi(D) = \dichi(D_{4n}') \ge n$. Indeed, note that any acyclic induced subdigraph $H'$ of $D$ satisfies $\omega(H') \leq 2$ by Lemma \ref{thm:triangle}.  Now, let $H$ be the underlying undirected graph of $H'$. Then the argument from the previous proof shows that $\chi(H) \le 4$. Since $\chi(G_{4n}) \ge 4n$ it follows that if $V(D)$ is partitioned into $t$ sets which induce acyclic subdigraphs, then $\chi(G_{4n}) \leq 4t$; therefore, we must have $t \ge n$. Thus, $\dichi(D) \ge n$, as claimed. This completes the proof.\bbox{}

\section{Further work}

Since $r=4$ is now the only open case of Conjecture \ref{conj:main}, it is natural to ask the following, which was first asked by James Davies (private communication):
\begin{question} 
Are there graphs with clique number 3 and arbitrarily large chromatic number whose all triangle-free induced subgraphs have chromatic number at most 3? 
\end{question}

Our construction was originally motivated by questions about the dichromatic number of graphs with bounded clique number. In view of Theorem \ref{thm:main2}, it is natural to ask: 

\begin{question}\label{q31}
Is there a function $f: \mathbb{N} \rightarrow \mathbb{N}$ such that for every digraph $D$ with no induced directed cycle of odd length, we have $\dichi(D) \leq f(\omega(D))$?
\end{question}

\begin{question}
\label{chordalquestion}
Is there a function $f: \mathbb{N} \rightarrow \mathbb{N}$ such that for every digraph $D$ with no induced directed cycle of length at least 4, we have $\dichi(D) \leq f(\omega(D))$?
\end{question}

 Question \ref{chordalquestion} asks about a directed analogue of chordal graphs and their dichromatic number. Since the first version of this paper, Question \ref{chordalquestion} has been answered in the negative by Aboulker, Bousquet, and de Verclos \cite{chordal}.
 We also ask the following more general question: 

\begin{question} \label{q33}
For which $l$ is there a function $f: \mathbb{N} \rightarrow \mathbb{N}$ such that for every digraph $D$ with no induced directed cycle of length not equal to $l$, we have $\dichi(D) \leq f(\omega(D))$?
\end{question}

Building on ideas of \cite{chordal}, we also answered Questions \ref{q31} and \ref{q33} in the negative \cite{secondpaper}.

\section*{Acknowledgments}

We are thankful to Louis Esperet for helpful comments on an earlier version of this paper, and thankful to the authors of \cite{chordal} for telling us about their result. 

We acknowledge the support of the Natural Sciences and Engineering Research Council of Canada (NSERC), [funding reference number RGPIN-2020-03912]. Cette recherche a été financée par le Conseil de recherches en sciences naturelles et en génie du Canada (CRSNG), [numéro de référence RGPIN-2020-03912]. Benjamin Moore is supported by the ERC-CZ project LL2005 (Algorithms and complexity within and
beyond bounded expansion) of the Ministry of Education of Czech Republic.

\end{document}